\documentclass[12pt]{article}
\usepackage{a4wide}
\usepackage{amssymb}
\usepackage{amsfonts}
\usepackage{amsmath}

\date{}

\newtheorem{proposition}{Proposition}[section]
\newtheorem{theorem}[proposition]{Theorem}
\newtheorem{lemma}[proposition]{Lemma}

\newtheorem{corollary}[proposition]{Corollary}

\def\der{\partial }

\def\nFM0{{\nu }_{F,M_0}}
\def\nFN0{{\nu }_{F,N_0}}
\def\nGN0{{\nu }_{G,N_0}}

\def\N0{ {\bf N}_0 }

\def\ra{\rightarrow}

\def\Xpm{X^{\pm }}

\def\l1{{\lambda}_1}

\def\a{\alpha}
\def\a0{ {\alpha }_0}
\def\a1{ {\alpha }_1}

\def\l{\lambda}

\def\nFGM0{{\nu }_{F,G,M_0}}

\def\nFN0{{\nu}_{F,N_0}}

\def\sm{{\sigma}^m}

\def\sm1{{\sigma}^{-1}}

\def\smtp1{{\sigma}^{-t+1}}

\def\S1{S^{-1}}

\def\Xpm1{X^{\pm 1}_1}

\def\sPM1{{\sigma }^{\pm 1}}
\def\sMP1{{\sigma }^{\mp 1 }}

\def\di{{\rm d.ind}}

\def\L{\Lambda}

\def\SL{\sum {}^n_{i=1}{\rm d}\Lambda {}_i}

\def\Ytm1{Y^{t-1}}
\def\Yim1{Y^{i-1}}

\def\h'{ \tilde{h} }

\def\Ass{{\rm Ass}}
\def\As{{\rm As}}
\def\ann{{\rm ann}}
\def\udim{{\rm u.dim }}
\def\ess{\subseteq_e}
\def\Spec{{\rm Spec}}
\def\lspec{{\rm l.spec}}
\def\rspec{{\rm r.spec}}
\def\mult{{\rm mult}}

\def\ufMod{{\rm ufMod}}
\def\Mod{{\rm Mod}}
\def\odo{\oplus\cdots\oplus}

\def\SL2Z{ {\rm SL}_2({\bf Z}) }

\def\Gp1{ G^{1 , 1 } }
\def\P11{ P^{-1 , 1 } }
\def\Pp1{ P^{1 , 1 } }

\def\nCLsr{{}^\nu\kern-2pt {\cal L}^{\sigma , \rho  }}
\def\nP{{}^\nu \kern-2pt P}
\def\nL{{}^\nu\kern-2pt L}
\def\nLL{{}^\nu\kern-2pt \Lambda}
\def\nPsr{{}^\nu\kern-2pt P^{\sigma , \rho  }}
\def\nLsr{{}^\nu\kern-2pt L^{\sigma , \rho  }}
\def\nuCL{{}^\nu\kern-2pt  {\cal L}}
\def\nCLsr{{}^\nu\kern-2pt {\cal L}^{\sigma , \rho  }}
\def\nCL1m{{}^\nu\kern-2pt {\cal L}^{-1 , 1  }}
\def\x1nu{x^\frac{1}{\nu}}
\def\xm1nu{x^{-\frac{1}{\nu}}}

\def\ra{\rightarrow }

\def\CI{{\cal I}}

\def\nAM0{{\nu }_{{\cal A},M_0}}
\def\nAN0{{\nu }_{{\cal A},N_0}}

\def\SL{{\rm SL}}

\def\Spec{{\rm Spec}}

\def\di!{\frac{\der^i}{i!}}
\def\dik!{\frac{\der^k_i}{k!}}

\def\N{\mathbb{N}}

\def\0{\overline{0}}
\def\1{\overline{1}}

\def\Ln1{\L_{n,\overline{1}}}

\def\a1{a_{\overline{1}}}

\def\S{\Sigma}

\def\vn1{\overrightarrow{n-1}}

\begin{document}

\author{V. V.\  Bavula 
}

\title{Primary Decompositions for Left Noetherian Rings}

\maketitle
\begin{abstract}
Two constructions are given that describe respectively all
shortest primary decompositions and all shortest uniform
decompositions for left Noetherian rings. They show that these
decompositions are, in general, highly non-unique.

{\em Key words: primary decomposition, uniform decomposition, the
left spectrum, left primes, uniform module.}

 {\em Mathematics subject classification
2000: 16P40, 16P70.}

$$ {\bf Contents} $$
\begin{enumerate}
\item Introduction. \item Primary decompositions for commutative
Noetherian rings. \item The left spectrum of a ring and associated
left primes. \item Primary decompositions. \item Description of
shortest primary decompositions. \item Description  of  shortest
uniform decompositions.

\end{enumerate}
\end{abstract}


\section{Introduction}

In this paper, module means a left module and all rings have 1.

Lasker-Noether's theory of the primary decomposition of a
submodule of a finitely generated  module over a commutative
Noetherian ring was generalized for modules over a not necessarily
commutative left Noetherian ring by L. Lesieur and R. Croisot
\cite{Les-Cro1963} and later by O. Goldman \cite{Goldman1969}. The
relations between Goldman's {\em primary} decomposition theory and
the {\em tertiary} decomposition theory of Lesieur and Croisot
were clarified by G. Michler \cite{Michler1970}. The former is a
`finer' decomposition theory than the later, it deals with modules
rather than with two-sided ideals. G. Michler showed that each
finitely generated Goldman-primary $R$-module is tertiary but not
vice versa, in general, and that each finitely generated tertiary
$R$-module is Goldman-primary iff non-isomorphic indecomposable
injective $R$-modules have different associated ideals.

In Section 2, we recall briefly main results on primary
decompositions in commutative situation. In Section 3, we collect
some results on the left spectrum of a ring and left primes. In
Section 4, it is proved that each submodule of uniformly finite
module admits a primary decomposition (Theorem \ref{RlufAss}).
Note that every finitely generated module over a left Noetherian
ring is uniformly finite.

 In Section 5, Theorem \ref{clminprdec} and Corollary \ref{almminprd}
describe respectively all the shortest primary decompositions and
all the maximal shortest primary decompositions of submodules of
uniformly finite modules.  They explain why primary decompositions
are not unique (in general, they are highly non-unique).

In Section 6,  (maximal shortest) uniform decompositions are
introduced. Theorem \ref{clstunid} and Corollary \ref{alsud}
describe respectively all the shortest uniform decompositions and
all the maximal  shortest uniform decompositions of submodules of
uniformly finite modules. It is proved that shortest uniform
decomposition is irredundant primary decomposition (Corollary
\ref{1sbudimMN}), and that each primary decomposition can be
refined to a uniform decomposition  (Lemma \ref{retpdu}).


\section{Primary decompositions for commutative Noetherian rings}

In this section,  $C$ is a {\bf commutative Noetherian ring} and
$M$ is a finitely generated $C$-module.

A prime ideal $P$ of $C$ is called  an {\bf associated prime} to
$M$ if $P$ is the annihilator of an element of $M$. The set of all
primes associated  to $M$ is written $\Ass_C(M)$. {\it $\Ass_C(M)$
is a finite nonempty set of primes each  containing the
annihilator $\ann_C(M)$ of the module $M$. The set $\Ass_C(M)$
includes all the primes minimal among primes containing
$\ann_C(M)$. The union of the associated primes of $M$ consists of
$0$ and the set of zerodivisors of $M$.} (\cite{EsenbudCA},
Theorem 3.1).

A submodule $N$ of  $M$ is {\bf primary} if $\Ass_C(M/N)$ consists
of just one prime, say $P$, then $N$ is called $P$-{\bf primary}.
The intersection of $P$-primary submodules is $P$-primary. {\bf
Primary decomposition} consists of writing an arbitrary submodule
$N$ of $M$ as a finite intersection of primary submodules.

The main properties of primary decomposition are given by the next
theorem (see \cite{EsenbudCA} for details). 

\begin{theorem}\label{comprdec}
Let $C$ be a commutative Noetherian ring and let $M$ be a finitely
generated $C$-module.
\begin{enumerate}
\item Any submodule $N$ of $M$ with $N\neq M$ is a finite
intersection of primary submodules, say $N=\cap_{i=1}^n N_i$ with
$N_i$ a $P_i$-primary submodule of $M$. \item Every associated
prime of $M/N$ occurs among the $P_i$. \item If the intersection
is {\bf irredundant} (i.e. no $N_i$ can be dropped) then the $P_i$
are precisely the associated primes of $M/N$. \item If the
intersection is {\bf shortest} (i.e. there is no such intersection
with fewer terms) then each associated prime of $M/N$ is equal to
$P_i$ for exactly one index $i$.
\end{enumerate}
\end{theorem}


\section{The left spectrum of ring and associated left primes}

In this section, the left prime spectrum of a ring and an analogue
of associates primes, so-called associated left primes (for
modules over noncommutative rings) are considered.

{\bf Essential submodules}. A submodule $N$ of a module $M$ is
{\bf essential} if $N\cap L\neq 0$ for all nonzero submodules $L$
of $M$, and we write $N\ess M$. The following (obvious) properties
of essential submodules are freely used in  proofs (\cite{MR},
2.2.2):
\begin{enumerate}
\item If $N\ess M$ and $M\ess L$ then $N\ess L$.
\item If $N\ess M$ and $L\ess M$ then $N\cap L\ess M$.
\item If $N_i\ess M_i$, $i\in I$, then $\oplus_{i\in I}N_i\ess \oplus_{i\in I}M_i$.
\item If $N$ is a submodule of $M$ then there exists a submodule
$N'$ of $M$ such that $N\cap N'=0$ and $N\oplus N'\ess M$.
\end{enumerate}

{\bf Uniform modules and  the uniform dimension}. A module $U$ is
{\bf uniform} if $U\neq 0$ and each nonzero submodule of $U$ is an
essential submodule. So, a nonzero module is uniform iff it does
not contain a direct sum of nonzero submodules. A module $M$ is
said  to have {\bf finite uniform dimension} ($\udim (M)< \infty
$) if it contains no infinite direct sum of nonzero submodules.
{\it If $M\neq 0$ has finite uniform  dimension then

$(i)$ $M$ contains an essential submodule which is a finite direct
sum, say $U_1\oplus \cdots \oplus U_n$, of uniform submodules
$U_i$ (\cite{MR}, 2.2.8).

$(ii)$ any direct sum of nonzero submodules of $M$ has at most $n$
summands, and

$(iii)$ a direct sum of uniform submodules of $M$ is essential in
$M$ iff it has precisely $n$ summands (\cite{MR}, 2.2.9).}

The nonnegative integer $n$ is called the {\bf uniform dimension}
(or {\bf Goldie dimension}) of $M$ and is written $\udim (M)$.

{\bf The essential equivalence $\sim_e $}. We say that modules
$M_1$ and $M_2$ are {\bf essentially equivalent} and write
$M_1\sim_e M_2$ if there exist essential submodules $L_1\ess M_1$
and $L_2\ess M_2$ with $L_1\simeq L_2$. It is obvious that this is
an equivalence relation. We denote by $[M]$ the equivalence class
of
 the module $M$ under $\sim_e$.

 {\it Definition}. The set of equivalence classes of
 all the uniform left  $R$-modules is called the {\bf left
 spectrum} $\lspec (R)$ of the ring $R$. Elements of the left
 spectrum are called {\bf left primes}.
 Similarly, the {\bf right
 spectrum} $\rspec (R)$ is defined.

Let $\widehat{R}$ be the set of isomorphism classes of simple left
$R$-modules. Clearly, two simple $R$-modules are isomorphic iff
they are essentially equivalent, hence $ \widehat{R}\subseteq
\lspec (R)$.

\begin{lemma}\label{l=r=SpR}
If $C$ is a commutative Noetherian ring then $\lspec (R)=\rspec
(R) =\{ [C/P]\, | $ $  \, P\in \Spec (C)\}$, and the map $ \lspec
(R)\ra \Spec (R)$, $[C/P]\mapsto P,$ is a bijection with the
inverse $P\mapsto [C/P]$.
\end{lemma}

{\it Proof}. This is evident. $\Box $

So, the left spectrum is a natural generalization of the spectrum
in the commutative situation.

For a left Noetherian ring $R$, let $\CI (R)$ be the set of
isomorphism classes of {\em indecomposable injective} $R$-modules.
The map $$\lspec (R)\ra \CI (R), \;\; [U]\mapsto [E(U)], $$ is a
bijection with the inverse $[E]\mapsto [E]$ where $E(U)$ is the
injective hull of the $R$-module $U$.

{\it Example}. Let $A_1=\langle x , \der \, | \, \der x -x\der =1
\rangle$ be the (first) {\em Weyl} algebra over a field $K$ of
characteristic zero. Note that the Weyl algebra $A_1$ is a simple
Noetherian domain which  contains no simple submodules, and any
proper factor module of the left $A_1$-module $A_1$ has finite
length. So, if $U$ is a uniform $A_1$-module then either $U$
contains a simple submodule, say $M$ (in this case, $[U]=[M]$) or
otherwise $U$ contains a copy
 of ${}_{A_1}A_1$ (in this case, $[U]=[A_1]$). It follows that
 $$ \lspec (A_1)=\widehat{A_1}\cup \{ [A_1]\}.$$

\begin{theorem}\label{multUV}
Let $M\neq 0$ be a module of finite uniform dimension $n$,
$\oplus_{i=1}^nU_i$ and $\oplus_{i=1}^nV_i$ be direct sums of
uniform submodules of $M$. Then, up to order of the $U_i$,
$[U_1]=[V_1], \ldots , [U_n]=[V_n]$.
\end{theorem}

{\it Proof}. Note that any direct sum of $n$ uniform submodules of
$M$ is an essential submodule of $M$. In particular, the
submodules $U:=\oplus_{i=1}^nU_i$ and $V:=\oplus_{i=1}^nV_i$ are
essential in $M$. We prove the following (more strong) statement:
{\em up to order of the $U_i$, there exist nonzero submodules
$V_1'\subseteq V_1, \ldots , V_n'\subseteq V_n$ such that
$[V_1']=[ U_1], \ldots , [V_n']=[U_n]$, and  for each} $i=1,\ldots
, n,$ $V_1'\oplus \cdots \oplus V_i'\oplus U_{i+1}\oplus \cdots
\oplus U_n\subseteq M$. Clearly, the intersection $V_1':=V_1\cap
U\neq 0$ since $V_1\neq 0$ and $U$ is an essential submodule of
$M$. For each $i=1,\ldots , n$, let $K_i$ be the kernel of the
canonical projection of $V_1'$ into $U_i$. Since $V_1'$ is a
submodule of $U$, the intersection $K_1\cap \cdots \cap K_n$ must
be zero, then one of the $K_i$ must be zero (otherwise, $K_1\cap
\cdots \cap K_n\neq 0 $ as $V_1'$ is a uniform module, a
contradiction). Up to order of the $U_j$, we may assume that
$K_1=0$. Then the projection map $V_1'\ra U_1$ is a monomorphism,
hence $[V_1']=[U_1]$ and
$$ V_1'+U_2+\cdots +U_n=V_1'\oplus U_2\oplus \cdots \oplus U_n.$$
This proves the claim  for  $i=1$.

We use induction on $i$. Suppose that $i>1$ and we have already
found required $V_1',\ldots , V_{i-1}'$. We have to find a nonzero
submodule $V_i'$ of $V_i$  such that $[V_i']=[U_i]$ and $
V_1'\oplus \cdots \oplus V_i'\oplus U_{i+1}\oplus \cdots\oplus
U_n\subseteq M$ (up to order of the $U_i, U_{i+1}, \ldots , U_n$).
Repeating the same argument as above for the submodule
$V_i':=V_i\cap (V_1'\oplus \cdots \oplus V_{i-1}'\oplus
U_{i}\oplus \cdots\oplus U_n)\neq 0$ of $V_1'\oplus \cdots \oplus
V_{i-1}'\oplus U_{i}\oplus \cdots\oplus U_n$ one can find $j\geq
i$ such that the projection $V_i'\ra U_j$ is a monomorphism (such
$j$ exists since $V_i'$ can not be a submodule of $V_1'\oplus
\cdots \oplus V_{i-1}'$). Up to re-numeration of the $U_i, \ldots
, U_n$ one can assume that $j=i$, hence $[V_i']=[U_i]$ and $
V_1'+\cdots +V_i'+U_{i+1}+\cdots +U_n=V_1'\oplus \cdots \oplus
V_i'\oplus U_{i+1} \oplus \cdots \oplus U_n$, as required. This
proves the claim and the theorem since $[V_1]=[V_1'],\ldots ,
[V_n]=[V_n']$. $\Box $

\begin{corollary}\label{1multUV}
 Let $M$ be a module of finite uniform dimension $n\geq 1$,
 $\oplus_{i=1}^n U_i$ be a direct sum of uniform submodules of
 $M$, $L$ be a uniform submodule of $M$. Then $[L]=[U_i]$ for some
 $i$.
\end{corollary}

{\it Proof}. Repeat the same argument as in the proof of the
Theorem \ref{multUV} in the case $i=1$ but for the nonzero module
$L':=L\cap \oplus_{i=1}^n U_i$.  $\Box $

Let $M$ be a nonzero $R$-module of finite uniform dimension $n$,
and let $\oplus_{i=1}^nU_i \ess M$ be a direct sum of uniform
submodules of $M$. A left prime $X\in \lspec (R)$ is called {\bf
associated} to $M$ if there exists a uniform submodule $U$ of $M$
 such that  $[U]=X$. Abusing language, we say that the uniform module $U$
 is {\bf associated} to $M$.

The set of associated left primes of $M$,
$$ \As_R(M):=\{ [U]\in \lspec (R)\, | \, U\; {\rm is \; a \;
uniform\; submodule\; of \;  } M\}$$ contains {\em finitely many}
elements since $\As_R(M)=\{ [U_i]\, | \, i=1, \ldots , n\}$
(Theorem \ref{multUV} and Corollary \ref{1multUV}). Let
$[U_1],\ldots , [U_s]$ be the (distinct) elements of the set
$\As_R(M)$. The {\bf multiplicity} of $[U_i]$ in $M$ is the number
of summands in $\oplus_{i=1}^n U_i$ that are essentially
equivalent to $U_i$, denoted $\mult_M([U_i])$ or simply $\mult
(U_i)$. Clearly, $[M]=[U_1^{m_1}\oplus \cdots \oplus U_s^{m_s}]$
with $m_i:= \mult_M([U_i])$.

Let $C$ be a commutative Noetherian ring and $M$ be a nonzero
finitely generated $C$-module.  By Theorem \ref{multUV} and
Corollary \ref{1multUV}, the map
$$ \Ass_C(M)\ra \As_C(M), \; P\mapsto [C/P],$$
 is a bijection with inverse $[U]\mapsto \ann_C(U)$. So,  $\Ass_C(M)$ and $\As_C(M)$ are
 essentially the same object. This is not the case if the ring $C$
 is noncommutative.

\begin{lemma}\label{Asdirs}
Let $R$ be a ring.
\begin{enumerate}
\item Let $M_1, \ldots , M_n$ be $R$-modules of finite uniform
dimension. Then $$\As_R(M_1\oplus \cdots \oplus
M_n)=\As_R(M_1)\cup \cdots \cup  \As_R(M_n).$$ \item Let $0\ra
M_1\ra M\ra M_2\ra 0$ be a short exact sequence of $R$-modules
each of them has finite uniform dimension. Then $\As_R(M)\subseteq
\As_R(M_1)\cup \As_R(M_2)$. \item Let $M_1, \ldots , M_n$ be
nonzero submodules of finite uniform dimension of a module $M$
such that $\As_R(M_i)\cap \As_R(M_j)=\emptyset $ for all $i\neq
j$. Then $M_1+\cdots +M_n=M_1\oplus \cdots \oplus M_n$.
\end{enumerate}
\end{lemma}

{\it Proof}. $1$. For each $i$, let  $\oplus_{j\in J_i}
U_{i,j}\ess M_i$ be an essential direct sum of uniform submodules
$U_{i,j}$ of $U_i$. Then $\As_R (M_i)=\{ [U_{i,j}]\, | \,  j\in
J_i\}$, and $\oplus_{i,j}U_{i,j}\ess \oplus_iM_i$ implies $\As_R
(\oplus_{i=1}^nM_i)=\{ [U_{i,j}]\, | \,  1\leq i\leq n, j\in
J_i\}=\cup_{i=1}^n\As_R (M_i)$.

$2$. If $U$ is a uniform submodule of $M$ then either $U\cap
M_1\neq 0$ (in that case, $[U]=[U\cap M_1]\in \As_R (M_1)$) or
otherwise $U\cap M_1=0$, in that case $U$ can be identified with a
uniform submodule of $M_2$, hence $[U]\in \As_R (M_2)$. Therefore,
$\As_R (M)\subseteq \As_R (M_1)\cup \As_R (M_2)$.

$3$. Suppose that $n=2$. Then $N:=M_1\cap M_2=0$ since otherwise,
that is, $N\neq 0$, we would have $\As_R(N)\neq \emptyset $ since
$\udim (N)\leq \udim (M_1)<\infty$ and, by the second statement,
$\As_R(N)\subseteq \As_R(M_1)\cap \As_R(M_2)=\emptyset$, a
contradiction. So, the result is true for $n=2$.

Let $n>2$. We use induction on $n$. By induction, the sum
$A:=M_1+\cdots +M_{n-1}$ is a direct sum. Then, by the first
statement, $\As_R(A)=\cup_{i=1}^{n-1}\As_R(M_i)$, and then
$\As_R(A)\cap \As_R(M_n)=\emptyset$. Hence, $A+M_n=A\oplus
M_n=M_1\oplus \cdots \oplus  M_n$. $\Box $


\section{Primary decompositions}

In this section, the notions of primary submodule and of primary
decomposition of a submodule will be generalized to the case of
submodules of uniformly finite modules. The main properties of
primary decompositions for commutative Noetherian rings, Theorem
\ref{comprdec}, still hold in this situation (Theorem
\ref{RlufAss} and Corollary \ref{1RlufAss}).

{\bf Uniformly finite rings and modules}.

{\it Definition}. Let $R$ be a ring, an $R$-module $M$ is called
{\bf a uniformly finite} $R$-module if all factor modules of $M$
have {\em finite} uniform dimension. Clearly, {\em each Noetherian
module is uniformly finite}. In particular, {\em all finitely
generated modules over an arbitrary left Noetherian ring are
uniformly finite}.

Given any submodules $N$ and $L$ of $M$ with $N\cap L=0$, there is
a submodule $C$ of $M$ that contains $L$ and is {\em maximal} with
respect to the property that $N\cap C=0$ (\cite{MR}, 2.2.3). This
is called a {\bf complement} to $N$ in $M$ that contains $L$.
Clearly, $N\oplus C\subseteq_e M$.

\begin{lemma}\label{ufm13}
Let $R$ be a ring, and $M_1,\ldots , M_n$, $M$ be $R$-modules.
\begin{enumerate}
\item If $0\ra M_1\ra M\ra M_2\ra 0$ is a short exact sequence of
$R$-modules then $M$ is a uniformly finite $R$-module  iff $M_1$
and $M_2$ are uniformly finite $R$-modules.
 \item $\oplus_{i=1}^n M_i$ is a uniformly finite $R$-module iff
 $M_1, \ldots , M_n$ are uniformly finite $R$-modules.
 \item If $M_1,\ldots , M_n$ are submodules of $M$ then
 $\sum_{i=1}^nM_i$ is a uniformly finite $R$-module iff
 $M_1, \ldots , M_n$ are uniformly finite $R$-modules.
\end{enumerate}
\end{lemma}

{\it Proof}. $1$. It is obvious that if $M$ is uniformly finite
then so are $M_1$ and $M_2$.

Suppose that the modules $M_1$ and $M_2$ are uniformly finite.
Then they have finite uniform dimension. Without loss of
generality we can assume that $M_2=M/M_1$. Given a submodule $N$
of $M$. We have to prove that $\udim (M/N)<\infty$. Consider the
short exact sequence of $R$-modules
$$ 0\ra M_1/M_1\cap N\ra M/N\ra M/(M_1+N)\simeq
(M/M_1)/(M_1+N/M_1)\ra 0.$$ The modules $M_1/M_1\cap N$ and $
M/(M_1+N)$ are uniformly finite as factor modules of the uniformly
finite  modules $M_1$ and $M_2=M/M_1$ respectively. So, without
loss of generality we can assume that $N=0$, and we have to prove
that $M$ has finite uniform dimension.

Choose any complement submodule to $M_1$ in $M$, say $C$. Then
$M_1\cap C=0$ and $M_1\oplus C\ess M$. The module $C$ can be
identified with its image in $M_2$ under the module epimorphism
$M\ra M_2=M/M_1$ since $M_1\cap C=0$. Then $\udim (C)\leq \udim
(M_2)<\infty$, hence $\udim (M_1\oplus C)=\udim (M_1)+\udim
(C)<\infty$, and so $\udim (M)=\udim (M_1\oplus C)<\infty$ since
$M_1\oplus C\ess M$, as required. This proves that $M$ is a
uniformly finite module.

$2$. This statement follows from the first statement.

$3$. If the sum $\sum_{i=1}^nM_i$ is a uniformly finite module
then so are the modules $M_1,\ldots , M_n$, by the first
statement. If $M_1,\ldots , M_n$ are uniformly finite modules then
so is their direct sum, by statement 2. The sum $\sum_{i=1}^nM_i$
is a uniformly finite module as an epimorphic image of the
uniformly finite module $\oplus_{i=1}^nM_i$. $\Box $

By Lemma \ref{ufm13}, if $M$ is a uniformly finite module then so
is {\em every  subfactor} of it.

{\it Definition}. A  ring $R$  is called a {\bf left uniformly
finite ring} iff ${}_RR$ is a uniformly finite module iff all
finitely generated  $R$-modules are uniformly finite (by Lemma
\ref{ufm13}, this definition makes sense). Every {\em left
Noetherian ring is a left uniformly finite ring}.

\begin{corollary}\label{ufrin12}

\begin{enumerate}
\item Each epimorphic image of a left uniformly finite ring is a
left uniformly finite ring. \item A direct product of rings
$R_1\times \cdots \times R_n$ is a left uniformly finite ring iff
$R_1, \ldots , R_n$ are left uniformly finite rings.
\end{enumerate}
\end{corollary}

{\it Proof}. This is a particular case of Lemma \ref{ufm13}.
$\Box $

The category $\ufMod (R)$ of uniformly finite $R$-modules is the
full subcategory of the category $\Mod (R)$ of $R$-modules. The
category $\ufMod (R)$ contains all the simple $R$-modules and is
closed under submodules, factor modules, extensions, and finite
direct sums. The category $\ufMod (R)$ contains all Noetherian
$R$-modules.

{\bf Primary submodules}.

{\it Definition}. A submodule $N$ of a module $M$ is {\bf primary}
if $\udim (M/N)<\infty$ and $\As (M/N)=\{ X\}$ (i.e. the set $\As
(M/N)$ consists of a single element), we say also that $N$ is
$X$-{\bf primary}.

For commutative rings, this definition coincides with the
classical one modulo the identification of prime ideals with their
correspondent uniform modules (Lemma \ref{l=r=SpR}). For a
commutative ring, a finite intersection of $P$-primary submodules
is a $P$-primary submodule. The same result is also true for
noncommutative rings.

\begin{lemma}\label{intprim}
Suppose that $N_1, \ldots , N_s$ are $X$-primary submodules of a
module $M$ then so is their intersection.
\end{lemma}

{\it Proof}. Via the module monomorphism
$$M/\cap_{i=1}^sN_i\ra
\oplus_{i=1}^sM/N_i, \; m+\cap_{i=1}^s N_i \mapsto (m+N_1, \ldots
, m+N_s),$$
 the module $M/\cap_{i=1}^sN_i$ can be viewed as
a submodule of the direct sum. Then,  $$\udim (M/\cap_{i=1}^s
N_i)\leq \udim (\oplus_{i=1}^s M/N_i)= \sum_{i=1}^s \udim
(M/N_i)<\infty$$ and,  by  Lemma \ref{Asdirs},
$$
 \As (M/\cap_{i=1}^sN_i)\subseteq  \As
(\oplus_{i=1}^sM/N_i)=\cup_{i=1}^s\As (M/N_i)=\{ X\},
$$
 so the
intersection $\cap_{i=1}^sN_i$  is an $X$-primary module.
 $\Box $

{\it Definition}. Let $N, N_1, \ldots , N_s$ be submodules of $M$
such that $N=\cap_{i=1}^s N_i$. The intersection $N=\cap_{i=1}^s
N_i$ is called a {\bf primary decomposition} of $N$ if each  $N_i$
 is $X_i$-primary submodule of $M$. The left primes $X_1, \ldots ,
 X_s$ are called the {\bf associated left primes} of the primary
 decomposition $N=\cap_{i=1}^s N_i$.

 In general, the set $\{ X_1, \ldots , X_s\}$ is not an invariant
 of the submodule $N$ of $M$ (it may vary for different primary
 decompositions), and some of the $X_i$'th may coincide. If the
  primary decomposition $N=  \cap_{i=1}^s N_i$ is {\bf irredundant}
 (i.e. no $N_i$ can be dropped) then the set $\{ X_1, \ldots ,
 X_s\}$ coincides  with the set $\As_R(M/N)$ of associated left
 primes of the module $M/N$ and the elements $X_1, \ldots , X_n$
  are {\em distinct} (Theorem \ref{RlufAss}), and so it is
 an invariant of the submodule $N$ of $M$. In this case, we say
 that $X_1, \ldots , X_n$ are the {\bf associated left primes} of
 the submodule $N$ of $M$, i.e. the associated left primes of the
 submodule $N$ of $M$ is the associated left primes of any {\em
 irredundant} primary decomposition of $N$ in $M$.

 For a finite set $S$, let $|S|$ be the number of elements in $S$.

\begin{theorem}\label{RlufAss}
Let $R$ be a ring, $M$ be a nonzero uniformly finite $R$-module.
 Then each  submodule $N$ of $M$ with $N\neq M$ is the
intersection of primary submodules, say $N=N_1\cap \cdots \cap
N_s$, with $\As_R(M/N_i)=\{ X_i\}$ (i.e. each submodule $N$ of $M$
such that $N\neq M$ admits a primary decomposition) and
\begin{enumerate}
\item $\As_R(M/N)\subseteq \cup_{i=1}^s\As_R(M/N_i)=\{ X_1, \ldots
, X_s \}$. \item If the intersection is irredundant then
$\As_R(M/N)=\cup_{i=1}^s\As_R(M/N_i)$. \item  If the intersection
is {\bf shortest} (i.e. there is no such intersection with fewer
terms) then each associated left prime of $M/N $  is equal to
$X_i$ for exactly one index $i$. \item Let $N= L_1\cap \cdots \cap
L_r$ be a primary decomposition of $N$ in $M$. Then it is shortest
iff $r=|\As_R(M/N)|$.
\end{enumerate}
\end{theorem}

{\it Proof}. The fact that each submodule $N$ of $M$ such that
$N\neq M$ admits a primary decomposition follows from Proposition
\ref{Csminprd}.

 Passing to the factor module $M/N$ we may assume that
$N=0$. Since $0=N=\cap_{i=1}^s N_i$, we have the $R$-module
monomorphism
$$ M\ra \oplus_{i=1}^s M/N_i, \; m\mapsto (m+N_1, \ldots , m+N_s).$$
So, we can identify the module $M$ with its image in the direct
sum. By the assumption, $N_i$ is a $X_i$-primary submodule of $N$,
so there is an essential submodule $U_{i,1}\oplus \cdots \oplus
U_{i, n_i}$ of $M/N_i$ where the $U_{i, j}$ are uniform modules
with $[U_{i,j}]=X_i$. Then $U:=\oplus_{i,j}U_{i,j}$ is the
essential submodule of $\oplus_{i=1}^sM/N_i$. If $V$ is a uniform
submodule of $M$ then $0\neq V':=V\cap U$ is a uniform submodule
of $U$, hence $[V]=[V']=[U_{i,j}]=X_i$ for some $i$ and $j$.
Therefore, $\As_R(M/N)\subseteq \cup_{i=1}^s\As_R(M/N_i)=\{ X_1,
\ldots , X_s \}$. This proves the first statement.

Suppose that the intersection $0=N=N_1\cap \cdots \cap N_s$ is
irredundant, that is, $N_i':=\cap_{j\neq i}N_j\neq 0$ for each
$i=1, \ldots , s$. Then $N_i'\cap N_i=0$, and so, for each $i$,
$N_i'$ can be viewed as a  nonzero submodule of $M/N_i$, and so
$\As_R(N_i')=\As_R(M/N_i)$. Now, $\cup_{i=1}^s
\As_R(M/N_i)=\cup_{i=1}^s \{ [ N_i'] \}\subseteq \As_R(M)$. The
inverse inclusion holds by the first statement. This proves the
second statement.

Suppose that the intersection $0=N=N_1\cap \cdots \cap N_s$ is
shortest. Then the elements $X_1,\ldots , X_s$ are distinct
(otherwise, say $X_i=X_j$ for some $i\neq j$, and so $N'=N_i\cap
N_j$ were a $X_i$-primary submodule of $M$ (Lemma \ref{intprim}),
$0=N=N'\cap (\cap_{k\neq i,j}N_k)$ were the intersection of
primary submodules which would contradict to the minimality of
$s$).  This proves the third statement.

The forth statement follows from the third. $\Box $

{\it Definition}. A submodule $N\subset M$ is {\bf irreducible} if
$N$ is not the intersection of two strictly larger submodules.

\begin{corollary}\label{1RlufAss}
Each Noetherian module $M$ and, in particular, every  finitely
generated module over a left Noetherian ring is a uniformly finite
module, and so Theorem \ref{RlufAss} holds.
\end{corollary}

The fact that {\em every submodule $N$ of a Noetherian module $M$
with $N\neq M$ admits a primary decomposition} can be proved
directly using the Noetherian  condition in a similar fashion as
in the commutative case.

{\it Proof}. Every submodule of $M$ which is not equal to $M$ is a
finite intersection of irreducible submodules. Otherwise, since
  $M$ is Noetherian, one could choose a submodule, say $L$, of
$M$ maximal among those submodules that do not share this
property. In particular, $L$ itself is not irreducible, so it is
the intersection of two strictly larger submodules, say  $L_1$ and
$L_2$. By the maximality of $L$, both modules $L_1$ and $L_2$ are
finite intersections of irreducible submodules, and it follows
that $L$ is too. This contradiction proves the claim and shows
that $N=N_1\cap \cdots \cap N_s$ with $N_i$ irreducible. In order
to finish the proof it suffices to show that {\em any irreducible
submodule $L\subseteq N$ is primary} (then it will follow that
$N=N_1\cap \cdots \cap N_s$ is a primary decomposition of $N$ in
$M$). Since $L$ is irreducible, we must have $\udim (M/L)=1$
(otherwise, $\udim (M/L)>1$, and we would find two uniform
submodules, say
 $\overline{V_1}$ and $\overline{V_2}$, of $M/L$ with
 $\overline{V_1}\cap \overline{V_2}=0$. Let $V_i$ be the
 pre-image of $\overline{V_i}$ under the epimorphism $M\ra M/L$,
 $m\mapsto m+L$. Then the submodules $V_1$ and $V_2$ would strictly
 contain the module $L$ and $V_1\cap V_2=L$, which would contradict the
 irreducibility of $L$).  So, $M/L$ is a uniform module, hence $L$
 is a $[M/L]$-primary submodule of $M$. $\Box $

If $R$ is  a commutative Noetherian ring  we get the classical
primary decomposition.


\section{Description of shortest primary decompositions}

In this section, maximal shortest primary decompositions are
introduced. Theorem \ref{clminprdec} and Corollary \ref{almminprd}
describe respectively all the shortest primary decompositions and
all the maximal shortest primary decompositions of submodules of
uniformly finite modules.  They explain why primary decompositions
are not unique (in general, they are highly non-unique). In the
commutative situation, another reason for non-uniqueness is given
in \cite{EsenbudCA}, 3.7.

{\bf Maximal shortest primary decompositions}. Let $R$ be a
 ring, $M$ be a nonzero uniformly finite $R$-module, $N$
be a submodule of $M$ such that $N\neq M$, $N=N_1\cap \cdots \cap
N_s$ be a shortest primary decomposition of $N$ in $M$. Then
$\As_R(M/N_1)=\{ X_1\}, \ldots , \As_R(M/N_s)=\{ X_s\}$ are
distinct sets, and $\As_R(M/N)=\{ X_1, \ldots , X_s\}$
(Proposition \ref{Csminprd} and Theorem \ref{RlufAss}).

{\it Definition}. The shortest primary decomposition $N=N_1\cap
\cdots \cap N_s$ is {\bf maximal} (with respect to inclusion) if
given another shortest primary decomposition $N=N_1'\cap \cdots
\cap N_s'$ with $\As_R(M/N_i')=\As_R(M/N_i)$ and $N_i\subseteq
N_i'$ for all $i=1, \ldots , s$ then $N_1=N_1', \ldots ,
N_s=N_s'$.

We will soon prove that for the submodule $N$ a shortest primary
decomposition always exits (Proposition \ref{Csminprd}), and each
shortest primary decomposition $N=N_1\cap \cdots \cap N_s$ is {\em
contained} in a maximal shortest primary decomposition $N=N_1'\cap
\cdots \cap N_s'$ (Theorem \ref{clminprdec}) which means that
$\As_R(M/N_i)=\As_R(M/N_i')$ and $N_i\subseteq N_i'$ for all
$i=1,\ldots , s$. The last statement is obvious if the module $M$
is Noetherian.

Recall that for  any submodules $N$ and $L$ of $M$ with $N\cap
L=0$, there is a submodule $C$ of $M$ that contains $L$ and is
{\em maximal} with respect to the property that $N\cap C=0$
(\cite{MR}, 2.2.3). This is called a {\em complement} to $N$ in
$M$ that contains $L$. Then $N\oplus C\subseteq_e M$. In general,
the complement $C$ is not unique (a typical example is when the
ring is a field then there exist many subspaces $C$ of $M$ with
the required properties). If $N\neq 0$ then $N\ess M/C$ via $M\ra
M/C$, $m\mapsto m+C$  (it is evident due to maximality of $C$).

\begin{proposition}\label{Csminprd}
Let $M$ be a nonzero uniformly finite $R$-module, $N$ be a
submodule of $M$ with $N\neq M$, $W_1, \ldots , W_s$ be nonzero
submodules of $M/N$ such that $W_1\oplus \cdots \oplus W_s\ess
M/N$ and the sets $\As_R(W_1)=\{ X_1\} , \ldots , \As_R(W_s)=\{
X_s\}$ are distinct. For each $i=1,\ldots , s$, let
$\overline{C}_i$ be a complement to $W_i$ in $M/N$ that contains
the modules $W_j$, $j\neq i$; and let
$C_i:=\pi^{-1}(\overline{C}_i)$ where $\pi :M\ra M/N$, $m\mapsto
m+N$. Then
\begin{enumerate}
\item each $W_i$ is an essential submodule of $M/C_i$, hence $C_i$
is an $X_i$-primary submodule of $M$. \item $N=C_1\cap \cdots \cap
C_s$ is
 a maximal shortest primary decomposition of $N$ in $M$.
\end{enumerate}
\end{proposition}

{\it Remarks}. $1$. The assumption that $W_1\oplus \cdots \oplus
W_s\ess M/N$ can be weakened to $W_1+ \cdots + W_s\ess M/N$ since
the later together with the assumption that the elements $X_1,
\ldots , X_s$ are distinct imply the former (Lemma
\ref{Asdirs}.(3)).

$2$. For each submodule $N$ of a uniformly finite module $M$ such
that $N\neq M$, one {\em can} find nonzero submodules $W_1, \ldots
, W_s$ of $M/N$ such that $W_1\oplus \cdots \oplus W_n\ess M/N$
and the sets $\As_R(W_1)=\{ X_1\}, \ldots , \As_R(W_n)=\{ X_n\}$
are {\em distinct}. Since $n:=\udim (M/N)<\infty$, take a direct
sum $U_1\odo U_n\ess M/N$ of uniform submodules of $M/N$. If
$\As_R(M/N)=\{ X_1, \ldots , X_s\}$ then, for  each $i=1,\ldots ,
s$, set $W_i:=\oplus \{  U_j\, | \, [U_j]=X_i\}$.

{\it Proof}. Passing to the factor module $M/N$, one can assume
that $N=0$. Then $\overline{C}_i= C_i$ is a complement to $W_i$ in
$M$, and so $W_i$ is an essential submodule of $M/C_i$ (via $M\ra
M/C_i$, $m\mapsto m+C_i$), hence $C_i$ is an $X_i$-primary
submodule of $M$ since $\As_R(W_i)=\{ X_i\}$. This proves the
first statement.

Let us show that $C:=C_1\cap \cdots \cap C_s=0$. Suppose that
$C\neq 0$, we seek a contradiction. Then $C':= C\cap (W_1\oplus
\cdots \oplus W_s)\neq 0$ since $W_1\oplus \cdots \oplus W_s$ is
 an essential submodule of $M$. Choose a nonzero element, say $c$, of $C'$. Then
$c=w_1+\cdots +w_s$ with $w_i\in W_i$. There exists $i$ such that
$w_i\neq 0$. Then $0\neq w_i=c-\sum_{j\neq i}w_i\in W_i\cap
C_i=0$, a contradiction. Therefore, $C=0$. Then the intersection
$C_1\cap \cdots \cap C_s=0$ is a primary decomposition of the zero
 module in $M$ (by the first statement). Now,
 the proof of Theorem \ref{RlufAss} is complete.
 The primary decomposition $C_1\cap \cdots \cap C_s=0$   is shortest (by
Theorem  \ref{RlufAss})  since $X_1, \ldots , X_s$ are distinct
and $\As_R(M)=\{ X_1, \ldots , X_s\}$.

Suppose that the shortest primary decomposition $0=C_1\cap \cdots
\cap C_s$ is contained in a shortest primary decomposition
$0=C_1'\cap \cdots \cap C_s'$, that is,
$\As_R(M/C_i)=\As_R(M/C_i')$ and $C_i\subseteq C_i'$ for all
$i=1,\ldots , s$. If there exists $i$ such that $C_i\neq C_i'$
then $C_i'\cap W_i\neq 0$, by the maximality of $C_i$, hence
$$0=C_1'\cap \cdots \cap C_s'\supseteq C_i'\cap (\cap_{j\neq i}C_j)\supseteq
C_i'\cap W_i\neq 0,$$
 a contradiction. Therefore, $0=C_1'\cap \cdots \cap C_s'$ is a maximal
 shortest primary decomposition of $0$ in $M$.  $\Box $

So, for any submodule of uniformly finite module there exists a
(maximal) shortest primary decomposition.

We will see shortly that Proposition \ref{Csminprd} gives {\em
all} the maximal shortest primary decompositions of  $N$ in $M$
(Corollary \ref{almminprd}).

\begin{lemma}\label{intirr}
Let $N_1, \ldots , N_s$ be nonzero submodules of a nonzero module
$M$ with $N_1\cap \cdots \cap N_s=0$. The intersection $N_1\cap
\cdots \cap N_s=0$ is irredundant iff $U_i:=\cap_{j\neq i}N_j\neq
0$ for all $i=1,\ldots , s$; in this case, $U_1+\cdots
+U_s=U_1\oplus \cdots \oplus U_s$.
\end{lemma}

{\it Proof}. The first statement is obvious, the second follows
from the fact that, for each $i$, $U_i\cap \sum_{j\neq
i}U_j\subseteq U_i\cap N_i=N_1\cap \cdots \cap N_s=0$. $\Box $

\begin{theorem}\label{clminprdec}
({\bf Description of shortest primary decompositions}) Let  $M$ be
a nonzero uniformly finite $R$-module, $N$ be a submodule of $M$
with $N\neq M$, $N=N_1\cap \cdots \cap N_s$ be a shortest primary
decomposition of $N$ in $M$,  $\As_R(M/N_i)=\{ X_i\}$. Then
\begin{enumerate}
\item there exist nonzero submodules $W_1, \ldots , W_s$ of $M/N$
such that $W_1\oplus \cdots \oplus W_s\ess M/N$ and $\As_R(W_i)=\{
X_i\}$ for each $i=1,\ldots , s$, and  \item for each $i=1,\ldots,
s$, there exists a complement, say $C_i$, of $W_i$ in $M/N$
containing $\sum_{j\neq i}W_j$ such that $\sum_{j\neq
i}W_j\subseteq N_i/N\subseteq C_i$.

In this way, all the shortest primary decompositions of $N$ in $M$
are obtained. In more detail,
 \item given any nonzero submodules $W_1',\ldots , W_s'$ of $M/N$
 such that $W_1'\oplus \cdots \oplus W_s'\ess M/N$
and $\As_R(W_i')=\{ X_i\}$ for each $i=1,\ldots , s$ (note that
$\As_R(M/N)=\{ X_1, \ldots , X_s\} )$, \item for each $i=1,\ldots,
s$, given a complement submodule $C_i'$ of $W_i'$ in $M/N$ that
contains the sum $\sum_{j\neq i}W_j'$, \item for each $i=1,\ldots
, s$, given a submodule $N_i'$ of $M$ containing $N$ and such that
$\As_R(M/N_i')=\{ X_i\}$ and  $\sum_{j\neq i}W_j'\subseteq
N_i'/N\subseteq C_i'$.
\end{enumerate}
Then $N=N_1'\cap \cdots \cap  N_s'$ is a shortest primary
decomposition of $N$ in $M$.
\end{theorem}

{\it Proof}. Passing to the factor module $M/N$, we can assume
that $N=0$. The shortest primary decomposition $0=N_1\cap \cdots
\cap N_s$ is irredundant. By Lemma \ref{intirr}, $W_i:=\cap_{j\neq
i}N_j\neq 0$ for all $i=1,\ldots , s$, and $W_1+\cdots
+W_s=W_1\oplus \cdots \oplus W_s$. For each $i$, $W_i\cap
N_i=N_1\cap \cdots \cap N_s=0$, hence $W_i\subseteq N/N_i$, and so
 $\As_R(W_i) = \{ X_i \}$. We claim that the sum $W:=W_1\oplus
\cdots \oplus W_s$ is an {\em essential} submodule of $M$. Suppose
that this is not true, then one can find a uniform submodule, say
$U$, of $M$ with $W\cap U=0$ (any nonzero module of finite uniform
dimension contains a uniform module). Since $\As_R(M)=\{
X_1,\ldots , X_s\}$ and $\As_R(W_1)=\{ X_1\}, \ldots ,
\As_R(W_s)=\{ X_s\}$, there exists $i$ such that $[U]=X_i$. Note
that, for each $j\neq i$, $U_j:=N_j\cap U\neq 0$ (since $N_j\cap
U=0$ would imply $U\subseteq M/N_j$, and then $\{ X_i\}
=\As_R(U)\subseteq \As_R(M/N_j)=\{ X_j\}$, this would contradict
to the fact that $X_i\neq X_j$). Each $U_j$ $(j\neq i)$ is a
nonzero submodule of the uniform module $U$, hence
$$0\neq \cap_{j\neq i}U_j\subseteq U\cap
 (\cap_{j\neq i}N_j)=U\cap W_i\subseteq U\cap W=0,$$
 a contradiction. This contradiction proves that $W$ is an
 essential submodule of $M$, and the first statement holds.

 Recall that $W_i\cap N_i=0$ for each $i$. Let $C_i$ be a
 complement to $W_i$ in $M$ such that $N_i\subseteq C_i$. In
 particular, the module $C_i$ contains the sum $\sum_{j\neq
 i}W_j\; (\subseteq N_i)$. So, the second statement holds.

It remains to show that all the shortest primary decompositions of
 $N=0$ in $M$ are obtained in this way. Let the modules
 $W_1',\ldots , W_s', C_1',\ldots , C_s', N_1',\ldots , N_s'$ be
 chosen as in statements 3-5. By Proposition \ref{Csminprd},
 $\cap_{i=1}^sN_i'\subseteq \cap_{i=1}^sC_i'=0$, and so
 $\cap_{i=1}^sN_i'=0$ is a shortest primary decomposition of $0$
 in $M$ since $\As_R(M/N_1')=\{ X_1\}, \ldots , \As_R(M/N_s')=\{
 X_s\}, \Ass_R(M)=\{ X_1, \ldots , X_s\}$, and $X_1, \ldots , X_s$
 are distinct (Theorem \ref{RlufAss}).
$\Box $

\begin{corollary}\label{almminprd}
({\bf Description of maximal shortest primary decompositions}) Let
$R$, $M$ and $N$ be as in Proposition \ref{Csminprd}. Then
Proposition \ref{Csminprd} describes all the maximal shortest
primary decompositions of $N$ in $M$.
\end{corollary}

{\it Proof}. This follows immediately from  Proposition
\ref{Csminprd} and Theorem \ref{clminprdec}. $\Box $

Let $R$ be a ring, $M$ be a nonzero uniformly finite $R$-module,
$N$ be a submodule of $M$ such that $N\neq M$, $\As_R(M/N)=\{ X_1,
\ldots , X_s\}$. The following algorithm explains {\em how one can
get a typical maximal shortest primary decomposition of $N$ in
$M$}.

\begin{enumerate}
\item Take any essential direct sum of uniform submodules in
$M/N$, say $U_1\oplus \cdots \oplus U_t\ess M/N$ (note that
$t=\udim (M/N)$). \item For each $i=1, \ldots , s$, let
$W_i:=\oplus \{ U_j\, | \, [U_j]=X_i\}$. Then $W_1\oplus \cdots
\oplus W_s\ess M/N$ and $\As_R(W_1)=\{ X_1\}, \ldots ,
\As_R(W_s)=\{ X_s\}$. \item For each $i=1, \ldots , s$, find a
complement, say $ \overline{C}_i$, to $W_i$ in $M/N$ that contains
$\oplus_{j\neq i}W_j$, and then take its preimage $C_i:=\pi^{-1}(
\overline{C}_i)$ under the module epimorphism $\pi :M\ra M/N$,
$m\mapsto m+N$. \item Then $N=C_1\cap \cdots \cap C_s$ is a
maximal shortest primary decomposition of $N$ in $M$, and all the
maximal shortest primary decompositions are obtained in this way.

The next step {\em gives a typical shortest primary decomposition
of $N$ in $M$}. \item For each $i=1, \ldots , s$, take a submodule
$N_i$ of $M$ containing $N$ and such that $\As_R(M/N_i)=\{ X_i\}$
and $\oplus_{j\neq i}W_j\subseteq N_i/N\subseteq C_i$ (eg.
$N_i=C_i$). Then $N=N_1\cap \cdots \cap N_s$ is a shortest primary
decomposition of $N$ in $M$ (and {\em all} shortest primary
decompositions of $N$ in $M$ are obtained in this way for all
possible choices of $N=C_1\cap \cdots \cap C_s$).
\end{enumerate}


\section{Description  of  shortest uniform  decompositions}

In this section, (maximal shortest) uniform decompositions are
introduced. Theorem \ref{clstunid} and Corollary \ref{alsud}
describe respectively all the shortest uniform decompositions and
all the maximal  shortest uniform decompositions of submodules of
uniformly finite modules. It will be proved that shortest uniform
decomposition is irredundant primary decomposition (Corollary
\ref{1sbudimMN}), and that each primary decomposition can be
refined to a uniform decomposition  (Lemma \ref{retpdu}).

In this section, {\em let $R$ be a ring, $M$ be a nonzero
uniformly finite $R$-module, $N$ be a submodule of $M$ such that
$N\neq M$, and $N=N_1\cap \cdots \cap N_s$ where $N_1, \ldots ,
N_s$ are submodules of $M$. }

{\bf Uniform decompositions}.

{\it Definition}. A decomposition $N=N_1\cap \cdots \cap N_s$ is
called {\bf a uniform decomposition} of $N$ in $M$ if each $M/N_i$
is a {\em uniform} $R$-module.

Obviously, each uniform decomposition is a primary decomposition
but,  in general, not vice versa.

\begin{lemma}\label{sbudimMN}
If $N=N_1\cap \cdots \cap N_s$ is a uniform decomposition of $N$
in $M$ then $s\geq \udim (M/N)$.
\end{lemma}

{\it Proof}. By the assumption, each $M/N_i$ is a uniform
$R$-module, so $\udim (M/N_i)=1$. Consider the $R$-module
monomorphism
$$ M/N\ra \oplus_{i=1}^sM/N_i, \;\; m+N\mapsto (m+N_1, \ldots ,
m+N_s).$$ Then $\udim (M/N)\leq \udim
(\oplus_{i=1}^sM/N_i)=\sum_{i=1}^s\udim (M/N_i)=s$. $\Box $

{\bf Shortest uniform decompositions}. A uniform decomposition
 $N=N_1\cap \cdots \cap N_s$ is {\bf shortest} if $s=\udim (M/N)$
 (the least possible value for $s$, by Lemma \ref{sbudimMN}).

\begin{corollary}\label{1sbudimMN}
A shortest uniform decomposition is an irredundant primary
decomposition.
\end{corollary}

{\it Proof}. This follows directly  from Lemma \ref{sbudimMN}.
$\Box $

 The next result proves existence of shortest uniform decompositions.

\begin{lemma}\label{UiCipr}
Let  $M$ be a nonzero uniformly finite $R$-module, $N$ be a
submodule of $M$ such that  $N\neq M$, $U_1, \ldots , U_n$ be
uniform submodules of $M/N$ such that $U_1\oplus \cdots \oplus
U_n$ is an essential submodule of $M/N$. For each $i=1, \ldots ,
n$, let $\overline{C}_i$ be a complement to $U_i$ in $M/N$ that
contains the modules $U_j$, $j\neq i$; and let
$C_i:=\pi^{-1}(\overline{C}_i)$ where $\pi :M\ra M/N$, $m\mapsto
m+N$. Then
\begin{enumerate}
\item each $U_i$ is an essential submodule of $M/C_i$, hence $C_i$
is a $[U_i]$-primary submodule of $M$. \item $\cap_{i=1}^nC_i=N$
is a shortest uniform decomposition of $N$ in $M$. In particular,
$\cap_{i=1}^nC_i=N$ is an irredundant primary decomposition of $N$
in $M$.
\end{enumerate}
\end{lemma}

{\it Proof}. $1$. By the choice of $\overline{C}_i$, $U_i$ is an
essential submodule of $(M/N)/\overline{C}_i\simeq M/C_i$.

$2$.  Without loss of generality one can assume that $N=0$.
 Suppose that the intersection $C:=C_1\cap \cdots \cap C_n$ is
nonzero, then $C':=C\cap (U_1\oplus \cdots \oplus U_n)\neq 0$
since $U_1\oplus \cdots \oplus U_n$ is an essential submodule of
$M$. Choose a nonzero element, say $c$, of $C'$. Then
$c=u_1+\cdots +u_n$ with $u_i\in U_i$. There exists $i$ such that
$u_i\neq 0$. Then $0\neq u_i=c-\sum_{j\neq i}u_i\in U_i\cap
C_i=0$, a contradiction. Therefore, $C=0$,  and $\cap_{i=1}^n
C_i=0$ is a shortest uniform decomposition of the zero module in
$M$ (by the first statement and since $\udim (M)=n$). By Corollary
\ref{1sbudimMN}, $\cap_{i=1}^n C_i=0$ is an irredundant primary
decomposition of $0$ in $M$.  $\Box $

The next theorem describes all shortest uniform decompositions for
submodules of nonzero uniformly finite modules.

\begin{theorem}\label{clstunid}
({\bf Description of shortest uniform decompositions}) Let $M$ be
a nonzero uniformly finite $R$-module, $N$ be a submodule of $M$
such that  $N\neq M$, $N=N_1\cap \cdots \cap N_n$ be a shortest
uniform decomposition of $N$ in $M$.
 Then
\begin{enumerate}
\item the module $M/N$ contains a direct sum $U_1\oplus \cdots
\oplus U_n$ of uniform submodules with $\As_R(M/N_i)=\{ [U_i]\}$
for each $i=1,\ldots, n$, and \item for each $i=1,\ldots, n$,
there exists a complement, say $C_i$, of $U_i$ in $M/N$ containing
$\sum_{j\neq i}U_j$ such that $\sum_{j\neq i}U_j\subseteq
N_i/N\subseteq C_i$.

In this way, all the shortest uniform decompositions of $N$ in $M$
are obtained. In more detail,
 \item given any direct sum $U_1'\oplus \cdots \oplus U_n'\subseteq  M/N$
 of uniform submodules (note  that $n=\udim (M/N)$),
 \item for each $i=1,\ldots,
n$, given a complement submodule $C_i'$ of $U_i'$ in $M/N$ that
contains the sum $\sum_{j\neq i}U_j'$, and
 \item for each
$i=1,\ldots , n$, given a submodule $N_i'$ of $M$ containing $N$
and such that $M/N_i'$ is a uniform module and $\sum_{j\neq
i}U_j'\subseteq N_i'/N\subseteq C_i'$.
\end{enumerate}
Then $N=N_1'\cap \cdots \cap N_n'$ is a shortest uniform
decomposition of $N$ in $M$. In particular, $N=N_1'\cap \cdots
\cap  N_n'$ is an irreducible primary decomposition of $N$ in $M$.
\end{theorem}

{\it Proof}. Passing to the factor module $M/N$, we can assume
that $N=0$. The shortest uniform decomposition $0=N_1\cap \cdots
\cap N_s$ is an irredundant primary decomposition (Corollary
\ref{1sbudimMN}).  By Lemma \ref{intirr}, $U_i:=\cap_{j\neq
i}N_j\neq 0$ for all $i=1,\ldots , n$, and $U_1+\cdots
+U_n=U_1\oplus \cdots \oplus U_n\subseteq M$. Since $n=\udim (M)$,
each $U_i$ must be a uniform submodule of $M$, and $U_1\oplus
\cdots \oplus U_n\ess M$.

For each $i=1,\ldots , n$, $U_i\cap N_i=N_1\cap \cdots \cap N_s=0$
implies $\{ [U_i]\}\subseteq \As_R(M/N_i)$, and so
$\As_R(M/N_i)=\{ [U_i]\}$ since the set $\As_R(M/N_i)$ consists of
a single element. So, the first condition holds.

 Let $C_i$ be a
 complement to $U_i$ in $M$ such that $N_i\subseteq C_i$
 (recall that $U_i\cap N_i=0$). In
 particular, the module $C_i$ contains the sum $\sum_{j\neq
 i}U_j\; (\subseteq N_i)$. So, the second condition holds.

It remains to show that all the shortest uniform decompositions of
 $N=0$ in $M$ can be  obtained in this way. Let the modules
 $U_1',\ldots , U_n', C_1',\ldots , C_n', N_1',\ldots , N_n'$ be
 chosen as in statements 3-5. By Proposition \ref{UiCipr},
 $\cap_{i=1}^nN_i'\subseteq \cap_{i=1}^nC_i'=0$, and so
 $\cap_{i=1}^nN_i'=0$ is a shortest uniform decomposition of $0$
 in $M$ since $n=\udim (M)$ and all the $M/N_i$ are uniform (by
 statement 5).
$\Box $

\begin{corollary}\label{1clstun}
 Let $M$ be a nonzero uniformly finite $R$-module, $N$ be a submodule of $M$
such that  $N\neq M$, $N=N_1\cap \cdots \cap N_n$ be a shortest
uniform decomposition of $N$ in $M$, $X_i:= [M/N_i]\in \lspec
(R)$. Up to order, let $X_1, \ldots , X_s$ be all the distinct
left primes among $X_1, \ldots , X_n$; and for each $i=1, \ldots ,
s$, let $n_i$ be the number of times the element $X_i$ occurs in
the series $X_1, \ldots , X_n$.
 Then $\As_R(M/N) = \{ X_1, \ldots , X_s\}$ and $ n_i=
 \mult_{M/N}(X_i)$ for all $i=1, \ldots , s$.
\end{corollary}

{\it Proof}. We may assume that statements 1 and 2 of Theorem
\ref{clstunid} hold. Since $n= \udim (M/N)$ the submodule
$U_1\oplus \cdots \oplus U_n$ of $M/N$ must be essential, and so
$\As_R(M/N) = \{ X_1, \ldots , X_s\}$ and $ n_i=
 \mult_{M/N}(X_i)$ for all $i=1, \ldots , s$.
$\Box $

{\bf Maximal shortest uniform decompositions}. A shortest uniform
decomposition $N=N_1\cap \cdots \cap N_n$ of $N$ in $M$ is {\bf
maximal} (with respect to inclusion) if given another shortest
uniform decomposition $N=N_1'\cap \cdots \cap N_n'$ of $N$ in $M$
 with $N_i\subseteq N_i'$ for all $i=1,\ldots , n$, then
 $N_1=N_1',\ldots , N_n=N_n'$.

Lemma \ref{UiCipr} and Theorem \ref{clstunid} show that each
shortest uniform
 decomposition $N=N_1\cap \cdots \cap N_n$ of $N$ in $M$ is {\em
 contained} in a maximal shortest uniform decomposition, say
 $N=N_1'\cap \cdots \cap N_n'$, which means that $N_i\subseteq
 N_i'$ for $i=1, \ldots , n$. In this case, by Theorem
 \ref{clstunid}, $\As_R(M/N_i)=\{ [\cap_{j\neq i}N_j]\}=\{
 [\cap_{j\neq i}N_j']\}=\As_R(M/N_i')$ as $\cap_{j\neq
 i}N_j\subseteq \cap_{j\neq i}N_j'$ is the inclusion of uniform
 modules.

\begin{corollary}\label{alsud}
({\bf Description of maximal shortest uniform decompositions}) Let
 $M$ and $N$ be as in Lemma \ref{UiCipr}. Then Lemma \ref{UiCipr}
describes all the maximal shortest uniform decompositions of $N$
in $M$.
\end{corollary}

{\it Proof}. This follows directly from Lemma \ref{UiCipr} and
Theorem \ref{clstunid}. $\Box $

Recall that $M$ be a nonzero uniformly finite $R$-module, $N$ be a
submodule of $M$ such that $N\neq M$. The following algorithm
explains {\em how one can get a typical maximal shortest uniform
decomposition of $N$ in $M$}.

\begin{enumerate}
\item Take any essential direct sum of uniform submodules in
$M/N$, say $U_1\oplus \cdots \oplus U_n\ess M/N$ (note that
$n=\udim (M/N)$). \item For each $i=1, \ldots , n$, find a
complement, say $ \overline{C}_i$, to $U_i$ in $M/N$ that contains
$\oplus_{j\neq i}U_j$, and then take its preimage $C_i:=\pi^{-1}(
\overline{C}_i)$ under the module epimorphism $\pi :M\ra M/N$,
$m\mapsto m+N$. \item Then $N=C_1\cap \cdots \cap C_n$ is a
maximal shortest uniform decomposition of $N$ in $M$, and all the
maximal shortest uniform decomposition are obtained in this way.

The next step {\em gives a typical shortest uniform decomposition
of $N$ in $M$}. \item For each $i=1, \ldots , n$, take a submodule
$N_i$ of $M$ containing $N$ and such that  $\oplus_{j\neq
i}U_j\subseteq N_i/N\subseteq \overline{C}_i$ and $M/N_i$ is
uniform  (eg. $N_i=C_i$). Then $N=N_1\cap \cdots \cap N_s$ is a
shortest uniform decomposition of $N$ in $M$ (and all shortest
uniform decompositions of $N$ in $M$ are obtained in this way for
all possible choices of $N=C_1\cap \cdots \cap C_n$).
\end{enumerate}

A primary decomposition $N=\cap_{i=1}^tN_i$ of $N$ in $M$ is
called {\bf a refinement} of a primary decomposition
$N=\cap_{j=1}^sM_j$ of $N$ in $M$ if the set $\{ 1, \ldots, t\}$
is a disjoint union of its non-empty subsets $I_1, \ldots , I_s$
such that, for each $j=1,\ldots , s$, $M_j=\cap_{k\in I_j}N_k$ and
$\As_R(M/N_k)=\As_R(M/M_j)$ for all $k\in I_j$.

\begin{lemma}\label{retpdu}
For any primary decomposition of a submodule of a uniformly finite
module there exists a refinement which is a uniform decomposition.
\end{lemma}

{\it Proof}. Let $N=N_1\cap \cdots \cap N_s$ be a primary
decomposition of a submodule $N$ of a uniformly finite module $M$
such that $N\neq M$. For each $i=1,\ldots , n$, let
$N_i=\cap_{j=1}^{n_i}N_{i,j}$, $n_i=\udim (M/N_i)$,  be a shortest
uniform decomposition of $N_i$ in $M$. Then $N=\cap_{i,j}N_{i,j}$
is a refinement of the primary decomposition $N=\cap_{i=1}^sN_i$,
and the refinement  is a uniform decomposition.
 $\Box $

Department of Pure Mathematics

University of Sheffield

Hicks Building

Sheffield S3 7RH

UK

email: v.bavula@sheffield.ac.uk

\end{document}